   \documentclass[11pt,letterpaper]{amsart}
\usepackage{amsmath, amssymb, amsfonts, graphicx, pgfplots, tikz, tikz-cd, stmaryrd}

\newtheorem{Theorem}{Theorem}[section]
\newtheorem{prop}[Theorem]{Proposition}

\def\beq#1#2\eeq{%
        \begin{equation}%
        \label{#1}%
            #2%
        \end{equation}%
   }

\usepackage{color}
\usepackage{graphicx}
\usepackage{epstopdf}

\title[Algebraic topology of the Lagrange inversion]{Algebraic topology of the Lagrange inversion}

\author{V.M. Buchstaber}\address{Steklov Mathematical Institute and Moscow State University, Russia}
\email{buchstab@mi-ras.ru}

\author{A.P. Veselov}
\address{Department of Mathematical Sciences,
Loughborough University, Loughborough LE11 3TU, UK}
\email{A.P.Veselov@lboro.ac.uk}

\begin{document}

\maketitle

\begin{abstract}
The Lagrange inversion formula for power series is one of the classical formulas from analysis and combinatorics.
A nice geometric interpretation of this formula in terms of the Stasheff polytopes was discovered by Loday.
We show that it also admits a natural topological interpretation in terms of the Chern numbers of the complex projective space.
The proof is based on our earlier work on the Chern-Dold character in complex cobordism theory and leads to a new derivation of the Lagrange inversion formula.
We provide a similar interpretation of the multiplicative inversion formula in terms of Chern numbers of the smooth theta divisors.
In this relation we introduce a new formal group defined by the Catalan numbers and explain the topological meaning of the corresponding Hirzebruch genus.
Finally, we discuss a related general problem of when all Chern numbers of an algebraic variety are divisible by its Euler characteristic.
  \end{abstract}


\section{Introduction}

Let  
$f(x)=x+\sum_{n\geq 1}a_n x^{n+1}$ be a formal power series and $g(y)=y+\sum_{n\geq 1}b_n y^{n+1}$ be its inverse under the substitution: $g(f(x))\equiv x.$
The Lagrange inversion problem \cite{Lagrange}  is to express the coefficients $b_n$ in terms of the coefficients $a_n.$

The celebrated {\it Lagrange inversion formula} \cite{Gessel, Stanley} can be written as
\beq{Lagrange}
(n+1) b_{n}=[x^n]\left(\frac{x}{f(x)}\right)^{n+1},
\eeq
where, following Stanley \cite{Stanley}, for any formal power series $h(x)$ we denote by $[x^n]h(x)$ the $n$-th coefficient of $h(x).$
We will provide a topological proof of this formula using the theory of characteristic classes \cite{MS}. We will see that this is closely related to the classical isomorphism 
$$
\tau(\mathbb CP^n)\oplus 1\cong \eta^{\oplus (n+1)},
$$
where $\eta$ is the dual to the tautological line bundle over $\mathbb CP^n$ (see \cite{MS}).

One can also compute $b_n=L_n(a_1,\dots, a_n)$ recursively with certain polynomials $L_n(a_1,\dots, a_n)$ with integer coefficients, which we call the {\it Lagrange inversion polynomials}:
$$
L_1=-a_1, \, L_2=-a_2+2a_1^2,\,\, L_3=-a_3+5a_1a_2-5a_1^3,
$$
$$
L_4=-a_4+6a_1a_3+3a_2^2-21a_1^2a_2+14a_1^4.
$$
  
  Their coefficients have a remarkable geometric interpretation  \cite{AA-2023, Loday} (see also \cite{BV-2025}) as the numbers (with signs) of the corresponding faces of the Stasheff polytopes (or, associahedra)
and an algebro-geometric interpretation as the numbers of strata of the Deligne-Mumford moduli space $M_{0,n}$ \cite{McMullen}.

The aim of this paper is to provide a topological interpretation of these coefficients as certain characteristic numbers of the complex projective spaces.

Let $M^n$ be a complex manifold and $\lambda=(\lambda_1, \dots, \lambda_k), \, i_1 \geq \dots \geq i_k$ be a partition of $n=|\lambda|:=\lambda_1+\dots+\lambda_k$. Let $c_{\lambda}(\nu M^n)$
be the characteristic Chern classes of the {\it normal} bundle $\nu M^n$ corresponding to the {\it monomial symmetric functions} $m_\lambda$ (see \cite{MS}).
The corresponding {\it Chern number} $c^\nu_{\lambda}(M^{n})$ is defined as the value of the cohomology class $c_\lambda(\nu M^{n})$ on the fundamental cycle $\langle M^{n} \rangle$:
$
c^\nu_{\lambda}(M^{n}):= (c_\lambda(\nu M^{n}), \langle M^{n} \rangle).
$

Consider the following generating function of these numbers
\beq{gennu}
C^\nu(M^n, t):=\sum_{\lambda: |\lambda|=n}c^{\nu}_\lambda(M^{n})t_\lambda, \quad t_\lambda=t_{\lambda_1}t_{\lambda_2}\dots t_{\lambda_k}.
\eeq

\begin{Theorem}
For the complex projective space $\mathbb CP^n$ the generating function of the monomial Chern numbers of the normal bundle $\nu(\mathbb CP^n)$ is
\beq{form}
C^\nu(\mathbb CP^n, t)=(n+1)L_n(t_1,\dots,t_n),
\eeq
where $L_n$ are the Lagrange inversion polynomials.
\end{Theorem}

Comparing this with the results of  \cite{AA-2023, Loday}, as a corollary we have a geometric interpretation of these characteristic numbers in terms of the combinatorics of Stasheff polytopes.
In particular, we have the following formula for the top Chern number of the normal bundle of $\mathbb CP^n$
$$
c_n(\nu \mathbb CP^n)=(-1)^n(n+1)C_n= (-1)^n {2n \choose n},
$$
where $C_{n}=\frac{1}{n+1}{ 2n \choose n}$ are the Catalan numbers.

We derive this theorem from the results of our paper \cite{BV-2024} about the Chern-Dold character in complex cobordism theory.

Let $\Theta^n$ be the smooth theta divisor of a general principally polarised abelian variety $A^{n+1}$ and $[\Theta^n]$
be its complex cobordism class, which  does not depend on the choice of $A^{n+1}$.
From the results of \cite{BV-2024} we have the following cobordism interpretation of the Lagrange inversion: 
\beq{form2}
[\mathbb CP^n]=(n+1)L_n(\tau_1,\dots,\tau_n), \quad \tau_k=\frac{[\Theta^k]}{(k+1)!},
\eeq
where $L_n$ are the Lagrange inversion polynomials.

It is interesting that the multiplicative inversion also admits topological interpretation. 

Let  
$F(x)=1+\sum_{n\geq 1}f_n x^{n}$ be a power series and $G(x)=1+\sum_{n\geq 1}g_n x^{n}$ be its multiplicative  inverse: $F(x)G(x)\equiv 1.$

The coefficients $g_n$ can be found recursively as $g_n=M_n(f_1,\dots, f_n)$ with certain polynomials $M_n(z_1,\dots, z_n)$ with integer coefficients, which 
we call the {\it multiplicative inversion polynomials}: 
$$M_1=-z_1,\, M_2=-z_2+z_1^2, \,
M_3=-z_3+2z_1z_2-z_1^3,$$
$$M_4=-z_4+2z_1z_3+z_2^2-3z_1^2z_2+z_1^4.
$$
They can be given as the {\it Hessenberg determinants} (see \cite{Insel}):
\beq{Hessen}
M_n(z_1,\dots, z_n)= (-1)^n
\begin{vmatrix}
z_1 & 1 & 0 & \hdots & \hdots & 0 \\
z_2 & z_1 & 1 & \hdots & \hdots & 0 \\
z_3 & z_2 & z_1 & 1 & \hdots & 0 \\
\vdots & \vdots & \vdots & \ddots & \ddots & 0 \\
z_{n-1} &  z_{n-2} & \hdots & \hdots & z_1 & 1 \\
z_n & z_{n-1} & \hdots & \hdots & z_2 & z_1
\end{vmatrix}.
\eeq



Let 
$$
C^\tau(M^n, t):=\sum_{\lambda: |\lambda|=n}c_\lambda(\tau M^{n})t_\lambda
$$
be the generating function of the monomial Chern numbers of the tangent bundle $\tau M^n$ of a complex manifold $M^n.$

\begin{Theorem}
The generating function of the monomial Chern numbers of the tangent bundle $\tau \Theta^n$ of the theta divisor is
\beq{form2}
C^\tau(\Theta^n, t)=(n+1)!M_n(t_1,\dots,t_n), 
\eeq
where $M_n$ are the multiplicative inversion polynomials.
\end{Theorem}

Note that the generating function of the monomial Chern numbers of the tangent bundle $\tau(\mathbb CP^n)$ of the projective space has the form
\beq{form3}
C^\tau(\mathbb CP^n, t)=(n+1)N_n(t_1,\dots,t_n), 
\eeq
where $N_n(t_1,\dots,t_n)$ is the generating function of the numbers of the noncrossing partitions of $n$ and can be expressed explicitly in terms of the partial ordinary Bell polynomials (see formulas (\ref{tanBell}), (\ref{Nn}) below): $N_1=t_1,$
$$
 N_2=t_2+t_1^2, \,\, N_3=t_3+3t_1t_2+t_1^3, \,\, N_4=t_4+4t_1t_3+2t_2^2+6t_1^2t_2+t_1^4, \dots.
$$

This motivates us to introduce the formal groups $F_\tau, F_\nu$ and Hirzebruch genera corresponding to two homomorphisms $\Phi_\tau, \Phi_\nu: \Omega_U \to \mathcal A$ defined by
$\Phi_\tau([M^n])=C^\tau(M^n,t),\quad \Phi_\nu([M^n])=C^\nu(M^n,t),$
where $\Omega_U$ is the complex cobordism ring and $\mathcal A=\mathbb Z[t_1,t_2,\dots]$ is the graded algebra of polynomials of infinitely many variables with degree $\deg t_k=k$. 

\begin{Theorem}
The exponential and logarithm of the formal groups $F_\tau$ and $F_\nu$ have the forms 
\beq{exp}
f_\tau(z)=z+\sum_{k=1}^\infty M_k(t) z^{k+1}, \quad f_\nu(z)=z+\sum_{k=1}^\infty t_k z^{k+1},
\eeq
\beq{log}
g_\tau(u)=u+\sum_{k=1}^\infty N_k(t) u^{k+1}, \quad g_\nu(u)=u+\sum_{k=1}^\infty L_k(t) u^{k+1}.
\eeq
\end{Theorem}


Specializing all $t_k=1$ in $F_\tau$ we get an interesting formal group $G_{\mathfrak C}$ with
the exponential and logarithm given by
\beq{tau*int}
f_{\mathfrak C}(z)=z-z^2, \quad g_{\mathfrak C}(u)=\sum_{k=0}^\infty C_k u^{k+1}=\frac{1-\sqrt{1-4u}}{2},
\eeq
where $C_k$ are the Catalan numbers, and the multiplication law
\beq{catprodint}
u*v=u+v-\frac{1}{2}(1-\sqrt{1-4u})(1-\sqrt{1-4v}).
\eeq
For the corresponding Hirzebruch genus $\Phi_{\mathfrak C}$ we have
\beq{chit}
\Phi_{\mathfrak C}(M^n)=(-1)^nc_n(\nu M^n),
\eeq
where $c_{n}(\nu M^n)$ is the usual Chern number of the {\it normal bundle} of a complex  manifold $M^n.$ 

Note that for any complex manifold $M^n$ and its tangent bundle $\tau M^n$ $c_{n}(\tau M^n)=\chi(M^n)$, which implies that for the specialized formal group $F_\nu$ the Hirzebruch genus $\Phi_\nu(M^n)=(-1)^n \chi(M^n)$  is dual to the Euler characteristic (see details in Section 3).

Finally we discuss a related problem of when all Chern numbers of an algebraic variety are divisible by its Euler characteristic.
As we have seen, $\mathbb CP^n$ and $\Theta^n$ are the examples of such varieties, but there are other examples shown in Section 4.

\section{Topology of Lagrange and multiplicative inversions}

 In complex cobordism theory a very important role is played by  the formal commutative group $F$ of  the geometric complex cobordisms introduced by Novikov \cite{Nov-67} 
\beq{F}
u*v=u+v+\sum_{i,j\geq 1}a_{ij}u^iv^j, \,\,\, a_{ij}\in \Omega_U^{2-2i-2j}.
\eeq
It can be defined by its logarithm $g_F(u)=u+\sum_{k=1}^\infty c_ku^{k+1},$
satisfying the relation $g_F(u*v)=g_F(u)+g_F(v),$ which can be given explicitly by Mishchenko formula
\beq{Mis}
g_F(u)=u+\sum_{n=1}^\infty[\mathbb CP^n]\frac{u^{n+1}}{n+1}.
\eeq
An explicit form of the exponential of this group $f_F(z),$ which is the inverse of $g_F(u)$, was found only recently in \cite{BV-2024}.

\begin{Theorem} [\cite{BV-2024}] The exponential of the formal group of geometric complex cobordisms can be given as
\beq{BV}
f_F(z)=z+\sum_{n=1}^\infty[\Theta^n]\frac{z^{n+1}}{(n+1)!},
\eeq
where $\Theta^n$ is smooth theta divisor of a principally polarised abelian variety.
\end{Theorem}

As a corollary we have the following representation of the cobordism class of the complex manifolds in terms of the theta divisors.

\begin{Theorem} [\cite{BV-2024}] 
The cobordism class of any complex manifold $M^{n}$ can be expressed as the sum is over all partitions $\lambda=(\lambda_1,\dots, \lambda_k)$ of $n$
\beq{decom}
[M^{n}]=\sum_{\lambda: |\lambda|=n}c^{\nu}_\lambda(M^{n})\frac{[\Theta^\lambda]}{(\lambda+1)!}, \quad \Theta^\lambda:=\Theta^{\lambda_1}\times\dots \times\Theta^{\lambda_k},
\eeq
where
$(\lambda+1)!:=(\lambda_1+1)!\dots(\lambda_k+1)!$ and $c^{\nu}_\lambda(M^{n})$ are the monomial Chern numbers of the normal bundle $\nu M^{n}.$ 
\end{Theorem}

Now we are ready to prove our results. 


To derive Theorem 1.1 we simply apply Theorem 2.2 to $M^n=\mathbb CP^n$ 
\beq{decomCP}
[\mathbb CP^{n}]=\sum_{\lambda: |\lambda|=n}c^{\nu}_\lambda(\mathbb CP^{n})\frac{[\Theta^\lambda]}{(\lambda+1)!}
\eeq
and compare this with formulas (\ref{Mis}),  (\ref{BV}).

Note that in this proof as well as in the proofs in \cite{BV-2024} the Lagrange inversion formula was not used.
Let us show that as a corollary we have now a topological proof of the Lagrange inversion formula (\ref{Lagrange}).

 Using the splitting principle we define the {\it universal monomial characteristic class} of a complex vector bundle $\xi$ on $M^n$
$$\mathcal C(\xi, t):=\prod_{j=1}^n (1+x_jt_1+x_j^2 t_2+\dots +x_j^n t_n),$$
where $x_1,\dots, x_n$ are the corresponding Chern roots \cite{MS} and $t=(t_1,\dots, t_n)$ are the parameters. 
It can be rewritten as
$$
\mathcal C(\xi, t):=\sum_{\lambda: |\lambda|\leq n}c_\lambda(\xi) t_\lambda, \quad t_\lambda=t_{\lambda_1}t_{\lambda_2}\dots t_{\lambda_k},
$$
where the characteristic class $c_\lambda(\xi)$ corresponds to the monomial symmetric polynomial $m_\lambda=x_1^{\lambda_1}\dots x_n^{\lambda_n}+\dots,$ 
and satisfies the property
$$
\mathcal C(\xi_1\oplus \xi_2,t)=\mathcal C(\xi_1,t)\mathcal C(\xi_2,t)
$$
for any two complex vector bundles on $M^n.$ In particular, for the tangent bundle $\tau M^n$ and normal bundle $\nu M^n$ we have
\beq{multip}
\mathcal C(\tau M^n,t)\mathcal C(\nu M^n,t)\equiv 1.
\eeq

To compute $\mathcal C(\tau M^n,t)$ for the projective space $M^n=\mathbb CP^n$ it is standard to use the well-known isomorphism
\beq{Isom}
\tau \mathbb CP^n\oplus 1 \cong \eta^{\oplus (n+1)},
\eeq
where $\eta$ is the dual to the tautological line bundle on $\mathbb CP^n$ (see \cite{MS}) with the first Chern class $c_1(\eta)=x$ generating $H^2(\mathbb CP^n, \mathbb Z)\cong \mathbb Z.$
We have 
\beq{tan}
\mathcal C (\tau \mathbb CP^n,t)=\mathcal C (\eta^{\oplus (n+1)},t)=(1+xt_1+x^2t_2+\dots+x^nt_n)^{n+1}, \, x^{n+1}=0.
\eeq
This implies that the monomial characteristic Chern class of the normal bundle is
\beq{nor}
\mathcal C (\nu \mathbb CP^n,t)=(1+xt_1+x^2t_2+\dots+x^nt_n)^{-(n+1)}, \quad x^{n+1}=0.
\eeq

The generating function $C^\nu (\mathbb CP^n,t)$ of the corresponding monomial Chern numbers is the coefficient at $x^n$ in the right hand side:
\beq{nu}
C^\nu (\mathbb CP^n,t)=[x^n] (1+xt_1+x^2t_2+\dots+x^nt_n)^{-(n+1)}=[x^n] \left(\frac{x}{f(x)}\right )^{n+1}
\eeq
where $f(x)=x+x^2t_1+x^3t_2+\dots+x^{n+1}t_{n}.$ Comparing this with formula (\ref{form}) from Theorem 1.1, we get the Lagrange inversion formula.

As a by-product we have the following interpretation of the generating function of the monomial Chern numbers of the tangent bundle of $\mathcal CP^n$:
\beq{tanBell}
C (\mathbb CP^n,t)=\hat B_{2n+1,n+1} (1,t_1,\dots,t_n),
\eeq
where $\hat B_{k,l}(z), k\geq l$ are the {\it partial ordinary Bell polynomials} \cite{Comtet}. This follows from formula (\ref{tan}) and the definition of the partial ordinary Bell polynomials
$$
(\sum_{m\geq 1} z_m x^m)^k=\sum_{n\geq k}\hat B_{n,k}(z)x^n.
$$

Here is the explicit form of the first few partial ordinary Bell polynomials (taken from Comtet \cite{Comtet}, page 309):
$$
\hat B_{3,2}=2 z_1z_2, \,\, \hat B_{5,3}=3 (z_1^2z_3+z_1z_2^2), \,\,  \hat B_{7,4}=4 (z_1^3z_4+3z_1^2z_2z_3+z_1z_2^3),
$$
$$
\hat B_{9,5}=5 (z_1^4z_5+4z_1^3z_2z_4+2z_1^3z_3^2+6z_1^2z_2^2z_3+z_1z_2^4).
$$

 The corresponding generating functions of tangent Chern numbers are
\beq{tanBell_n}
C^\tau (\mathbb CP^n)=[x^n] (1+xt_1+x^2t_2+\dots+x^nt_n)^{n+1},
\eeq
or, more explicitly,
$
C^\tau (\mathbb CP^1)=2t_1, \,\, C^\tau (\mathbb CP^2)=3 (t_2+t_1^2), 
$
$$
\,\,  C^\tau (\mathbb CP^3)=4 (t_3+3t_1t_2+t_1^3), \,\, C^\tau (\mathbb CP^4)=5 (t_4+4t_1t_3+2t_2^2+6t_1^2t_2+t_1^4).
$$
In combinatorics the quotients 
\beq{Nn}
N_n(t_1,\dots,t_n):=\frac{C^\tau (\mathbb CP^n)}{n+1}=\frac{1}{n+1}[x^n] (1+xt_1+x^2t_2+\dots+x^nt_n)^{n+1}
\eeq 
are known as the generating functions of the numbers of noncrossing partitions\footnote{We are grateful to Tom Copeland for pointing this out to us.} (see Ex.6.19(pp) in Stanley \cite{Stanley} and the entry A134264 in the On-Line Encyclopedia of Integer Sequences (OEIS)). 
Namely, the coefficient of $N_n(t_1,\dots,t_n)$ at $t_\lambda$ gives the number of such partitions of the type $\lambda, \, |\lambda|=n$. For example, $N_3=t_3+3t_1t_2+t_1^3$ describes the numbers of noncrossing partitions of 3: $\{123\}, \{12-3, 13-2, 23-1\}, \{1-2-3\}.$

The formula (\ref{Nn}) provides the topological interpretation of these numbers. 
Note that when we specialize here all $t_i=1$ we get the Catalan number
$$N_n(1,\dots,1)=\frac{1}{n+1}[x^n] (1-x)^{-(n+1)}=\frac{1}{n+1} {2n \choose n}=C_n,$$
which, in particular, gives the number of vertices of $(n-1)$-dimensional Stasheff polytope (or, different triangulations of $(n+2)$-gon).

For the normal bundle of $\mathbb CP^n$ due to Theorem 1.1 the generating functions of Chern numbers are
$$
C^\nu (\mathbb CP^1)=-2t_1, \,\, C^\nu (\mathbb CP^2)=3 (-t_2+2t_1^2), \,\,  C^\nu (\mathbb CP^3)=4 (-t_3+5t_1t_2-5t_1^3),
$$
$$
C^\nu (\mathbb CP^4)=5 (-t_4+6t_1t_3+3t_2^2-21t_1^2t_2+14t_1^4).
$$
The last formula is an agreement with the fact that $3D$ Stasheff polytope has 6 pentagonal and 3 quadrilateral faces, 21 edges and 14 vertices (see Fig. 1 and \cite{AA-2023, BV-2025}).
\begin{figure}[h]
\includegraphics[width=40mm]{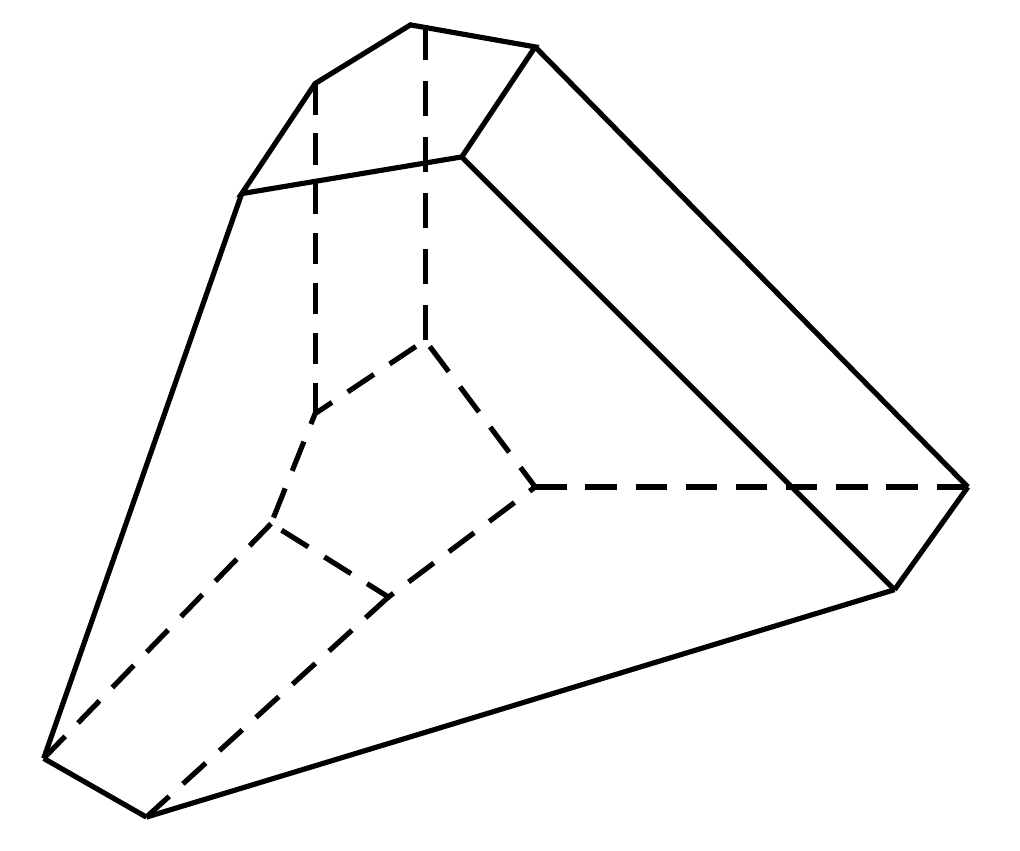} 
\caption{Stasheff polyhedron}
\end{figure} 

To prove Theorem 1.2 we note that the normal bundle $\nu\Theta^n$ can be identified with the line bundle $\mathcal L=i^*(L)$, where $i: \Theta^n \to A^{n+1}$ is natural embedding and $L$ is the principal polarisation line bundle on $A^{n+1}$ (see \cite{BV-2024}). 
This means that for the normal bundle the corresponding monomial Chern class is
\beq{normtheta}
\mathcal C(\nu \Theta^n,t)=1+\sum_{k=1}^n x^n t_n, \,\,\, x=c_1(\mathcal L).
\eeq
Now Theorem 1.2 follows from the relation
$\mathcal C(\tau \Theta^n,t)\mathcal C(\nu \Theta^n,t)\equiv 1$
and the fact that 
$
(x^n, \langle \Theta^{n} \rangle)=(n+1)!
$
(see \cite{BV-2024}). 

Thus the generating functions of the monomial Chern numbers of the tangent and normal bundles of $\Theta^n$ are
\beq{numtheta}
C^\nu(\Theta^n,t)=(n+1)!t_n, \quad
C^\tau(\Theta^n,t)=(n+1)!M_n(t_1,\dots,t_n), 
\eeq
in particular, $C^\tau (\Theta^1)=-2t_1, \,\, C^\tau (\Theta^2)=6 (-t_2+t_1^2),$
$$
C^\tau (\Theta^3)=24 (-t_3+2t_1t_2-t_1^3), \,\,\,
C^\tau(\Theta^4)=120 (-t_4+2t_1t_3+t_2^2-3t_1^2t_2+t_1^4).
$$

To link the corresponding formulas with combinatorics of permutohedra one should consider the multiplicative inversion for the series in the {\it exponential form}
$$P(x)=1+\sum_{n=1}^\infty p_{n} \frac{x^{n}}{n!}.$$
 The coefficients of the multiplicative inverse $$Q(x)=\frac{1}{P(x)}=1+\sum_{n=1}^\infty q_{n} \frac{x^{n}}{n!}$$ can be expressed as 
$
q_n=\hat M_n(p_1,\dots, p_n), 
$
where the polynomials $\hat M_n$ are simply related to the multiplicative inversion polynomials by the formula
\beq{hat}
\hat M_n(p_1,\dots, p_n)=n! M_n(z_1,\dots, z_n), \,\, z_k=\frac{p_k}{k!},
\eeq
$$
\hat M_1=-p_1,\, \hat M_2=-p_2+2p_1^2,
\hat M_3=-p_3+6p_1p_2-6p_1^3,
$$
$$ 
\hat M_4=-p_4+8p_1p_3+6p_2^2-36p_1^2p_2+24p_1^4.
$$
Their coefficients have a natural interpretation as the numbers of the corresponding faces of permutohedra (see \cite{AA-2023}).
In particular, the formula for $\hat M_4$ is in agreement with the fact that 3D permutohedron has 8 hexagonal and 6 quadrilateral faces, 36 edges and 24 vertices (see Fig. 2). 

\begin{figure}[h]
\includegraphics[width=35mm]{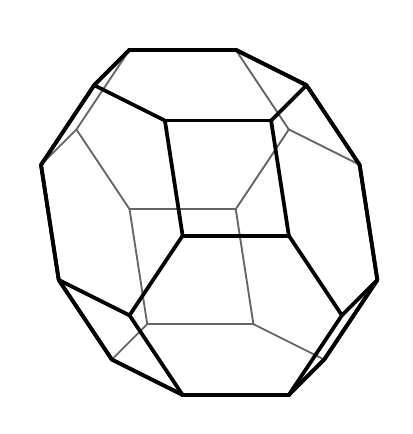} 
\caption{3D permutohedron.}
\end{figure}

\section{Related formal groups and Hirzebruch genera}

Recall that the Hirzebruch genus \cite{Hirz-66} in complex cobordisms \cite{B-1970} is a homomorphism $\Phi: \Omega_U \to \mathcal A$, where $\mathcal A$ is some algebra over $\mathbb Q$, determined by its characteristic power series
$
Q(z)=1+\sum_{n=1}^\infty a_n z^n, \quad a_n \in \mathcal A.
$

In \cite{BV-2024} we have proved that the exponential generating function of $\Phi(\Theta^n)$ is related to the corresponding $Q(z)$ by the formula
\beq{hirzn}
\frac{z}{Q(z)}=z+\sum_{n=1}^\infty \Phi(\Theta^n)\frac{z^{n+1}}{(n+1)!}.
\eeq

Let $\mathcal A=\mathbb Z[t_1,t_2,\dots,t_n,\dots]$ be the graded algebra of polynomials of infinitely many variables with degree $\deg t_k=k$. 

It is known that the complex cobordism ring $\Omega_U$ is multiplicatively generated by some complex manifolds.
Consider two natural homomorphisms $\Phi_\tau, \Phi_\nu: \Omega_U \to \mathcal A$ sending the cobordism class of a complex manifold $M^n$ to the generating function of the corresponding monomial Chern numbers of the tangent and normal bundles respectively:
\beq{hirztn}
\Phi_\tau([M^n])=C^\tau(M^n,t),\quad \Phi_\nu([M^n])=C^\nu(M^n,t) \in \mathcal A.
\eeq

\begin{prop}
The characteristic series of the Hirzebruch genera $\Phi_\tau$ and $ \Phi_\nu$ are respectively
\beq{chartn}
Q_\tau(z)=1+\sum_{k=1}^\infty t_k z^k, \quad Q_\nu(z)=1+\sum_{k=1}^\infty M_k(t) z^k,
\eeq
where $M_k(t_1,\dots,t_k)$ are the multiplicative inversion polynomials.
\end{prop}

Indeed, this immediately follows from formulas (\ref{numtheta}) and (\ref{hirzn}).

Consider now the formal group $F$ of geometric complex cobordisms (\ref{F}). 
For any homomorphism $\Phi: \Omega_U \to \mathcal A$ we can define a formal group $F_\Phi$ with coefficients in $\mathcal A.$
According to Quillen \cite{Quil-69} any formal group over $\mathcal A$ can be described this way.

In particular, for $\mathcal A=\mathbb Z[t_1,t_2,\dots,t_n,\dots]$ and $\Phi=\Phi_\tau, \Phi_\nu$ we have the corresponding formal groups $F_\tau$ and $F_\nu$ over $\mathcal A.$
Theorem 1.3 claims that the exponential and logarithm of the formal groups $F_\tau$ and $F_\nu$ have the form (\ref{exp}), (\ref{log}).
To prove this we can simply apply the homomorphisms $\Phi_\tau, \Phi_\nu$ to (\ref{Mis}), (\ref{BV}) and use Theorems 1.1, 1.2 and formulas (\ref{form3}), (\ref{numtheta}).

Let us consider now the specializations $F_{\tau, spec}$ and $F_{\nu, spec}$ of these groups by setting all $t_k=1.$ 

\begin{prop}
The exponential and logarithm of the specialised formal groups $F_{\tau, spec}$ and $F_{\nu, spec}$ are respectively
\beq{tau*}
f_{\tau, spec}(z)=z-z^2, \quad g_{\tau, spec}(u)=\sum_{k=0}^\infty C_k u^{k+1}=\frac{1-\sqrt{1-4u}}{2},
\eeq
\beq{nu*}
f_{\nu, spec}(z)=\sum_{k=0}^\infty z^{k+1}=\frac{z}{1-z}, \quad g_{\nu, spec}(u)=\sum_{k=0}^\infty (-1)^k u^{k+1}=\frac{u}{1+u},
\eeq
where $C_k$ are the Catalan numbers. 

The corresponding Hirzebruch genera are respectively
\beq{taunugen}
\Phi_{\tau,spec}(M^n)=(-1)^nc_n(\nu M^n), \quad \Phi_{\nu,spec}([M^n])=(-1)^n \chi(M^n),
\eeq
where $\chi(M^n)$ is the Euler characteristic and $c_{n}(\nu M^n)$ is the usual Chern number of the normal bundle of $M^n.$
\end{prop}

Indeed, when all $t_k=1$ we have the geometric series $F=\sum_{k=0}^\infty x^k=\frac{1}{1-x}$ with the multiplicative inverse $G=1-x$. This means that $M_k(1,\dots,1)=0$ for $k>1$ with $M_1(1)=-1,$ implying that $f_{\tau, spec}(z)=z-z^2.$ To find the inversion we simply solve the equation $x-x^2=u$ giving
$x=\frac{1-\sqrt{1-4u}}{2}$, which is well-known to be the generating function of the Catalan numbers (this, of course, also follows from the combinatorial interpretation of $N_k(t)$ with $N_k(1,\dots,1)=C_k$). This implies formula (\ref{tau*}), formula (\ref{nu*}) follows immediately from the definition.
 
To discuss the corresponding Hirzebruch genera we need the following general notion of duality in complex cobordisms.

Let $\Psi^{(-1)}: \Omega_U\to\Omega_U$ be the Adams-Novikov operation in complex cobordisms, which is the involution induced by complex conjugation in the universal bundle over classifying space $BU(n)$. It acts on each cobordism class of a complex manifold by the formula
$$\Psi^{(-1)}([M^n])=(-1)^n[M^n].$$
For any Hirzebruch genus $\Phi: \Omega_U \to \mathcal A$ the {\it dual} genus is defined as the composition $\Phi^*=\Phi\circ \Psi^{(-1)}$, so for complex manifolds
\beq{dual}
\Phi^*([M^n])=(-1)^n\Phi([M^n]).
\eeq
 The exponential and logarithm of the corresponding formal group and its dual are related by
\beq{dual}
f^*(z)=-f(-z), \quad g^*(z)=-g(-z).
\eeq
It is known that for the complex manifolds the Euler characteristic $\chi$ coincides with the Hirzebruch genus with $Q_\chi(z)=1+z$ (see e.g. \cite{BN-1971}). For the corresponding formal group $F_\chi$ we have 
$$
f_\chi(z)=\frac{z}{1+z}=f^*_{\nu,spec}(z), \quad g_\chi(u)=\frac{u}{1-u}=g^*_{\nu,spec}(u),
$$
so we see that $\Phi_{\nu, spec}=\chi^*$ is dual to the Euler characteristic and therefore
$$
\Phi_{\nu,spec}([M^n])=(-1)^n \chi(M^n).
$$

The case of the formal group $\mathfrak C=F_{\tau, spec}$, which we call the {\it Catalan formal group,} is more interesting. Its multiplication law has the form
$$u*v=u+v-\frac{1}{2}(1-\sqrt{1-4u})(1-\sqrt{1-4v})=u+v-2\sum_{k,l=0}^\infty C_kC_l u^{k+1}v^{l+1}.$$

For the corresponding {\it Catalan genus} $\Phi_{\mathfrak C}$ we have
$
\Phi_{\mathfrak C}(\mathbb CP^n)=(n+1) C_n, 
$
and 
$\Phi_{\mathfrak C}(\Theta^n)=0$ for $n>1$ with $\Phi_{\mathfrak C}(\Theta^1)=-2.$
Applying this to both sides of (\ref{decom}) we find that for any complex manifold $M^n$
\beq{chit}
\Phi_{\mathfrak C}(M^n)=(-1)^nc_n(\nu M^n),
\eeq
where $c_{n}(\nu M^n)$ is the usual Chern number of the normal bundle of $M^n.$ 

The dual Catalan genus 
\beq{chitdual}
\Phi^*_{\mathfrak C}(M^n)=c_n(\nu M^n),
\eeq
is a direct analogue of the Euler characteristic $\chi([M^n])=c_{n}(\tau M^n)$ with the tangent bundle of $M^n$ replaced by its normal bundle.

The exponential and logarithm of the corresponding dual Catalan formal group $\mathfrak C^*$ are respectively
\beq{tauC*}
f^*_{\mathfrak C}(z)=z+z^2, \quad g^*_{\mathfrak C}(u)=\sum_{k=0}^\infty (-1)^{k}C_k u^{k+1}=\frac{\sqrt{1+4u}-1}{2},
\eeq
the multiplication law is
$$u*v=u+v+\frac{1}{2}(\sqrt{1+4u}-1)(\sqrt{1+4v}-1)=u+v+2\sum_{k,l=0}^\infty C_kC_l u^{k+1}v^{l+1}.$$

\section{Divisibility of Chern numbers}

As we have seen from Theorem 1.1 all the Chern numbers of the complex projective space $\mathbb CP^n$ are divisible by $n+1$. This fact was observed by the first author already in 1970, who found it a bit surprising, taking into account that the Todd genus of $\mathbb CP^n$ is 1.
Note that the reduced generating function $C^\nu(\mathbb CP^n,t)/(n+1)=L_n(t_1,\dots, t_n)$ is precisely the corresponding Lagrange inversion polynomial.

This divisibility is not obvious from formula (\ref{tanBell}), but this is true since the greatest common divisor (gcd) of the coefficients of the partial ordinary Bell polynomial $\hat B_{k,l}(z)$ is known to be $l/gcd(k,l)$ (see \cite{O'Sullivan}). For $\hat B_{2n+1,n+1}$ we have $gcd(2n+1,n+1)=1,$ so it is also divisible by $n+1$.

A similar divisibility happens for the theta divisor $\Theta^n$, where all Chern numbers are divisible by $(n+1)!.$

A natural question is how often we have the divisibility of all Chern numbers. Since $c_n(M^n)=\chi(M^n)$ is the Euler characteristic for all complex manifolds, the greatest common factor $d$ of the Chern numbers $c_\lambda(M^n)$ must be a divisor of $\chi(M^n).$

We can view this as a part of the general Milnor-Hirzebruch problem, which was first posed in \cite{Hirz}. Its algebraic version can be formulated as follows (see \cite{BV-2024}):

{\it Characterise all possible sets of the integers, which can be realised as the Chern numbers $c_\lambda(M^n)$ of some irreducible complex algebraic variety $M^n$.}

In this version the problem is largely open, although some arithmetic restrictions of are known since the work of Milnor and Hirzebruch \cite{Hirz}.
The restriction to the irreducible varieties is essential, since without this assumption the problem was resolved (see \cite{Stong-68}). 

For any algebraic variety $M^n$ we define the number
$d(M^n)$ as the greatest common divisor of the corresponding Chern numbers $c_\lambda(M^n), \, |\lambda|=n.$ Equivalently, we can use here the Chern numbers of the normal bundle $c^\nu_\lambda(M^n)$ since they are related to the usual Chern numbers by an invertible linear change with integer coefficients.

The results of the previous section motivate the following problem:

{\it Characterise all irreducible algebraic varieties $M^n$ with $\chi(M^n)\neq 0$ such that $d(M^n)=|\chi(M^n)|.$}

 We call such varieties {\it numerically extremely divisible}.
We have seen that this is the case for $\mathbb CP^n$ and $\Theta^n$ with Euler characteristics $\chi(\mathbb CP^n)=n+1$ and $\chi(\Theta^n)=(n+1)!$ respectively.
We can extend the second series of examples as follows. 

Let $\Theta^n(k), k\geq 2$ be a smooth zero locus of generic section of the $k$-th tensor power $L^{\otimes k}$ of the line bundle $L,$ defining the principal polarisation of an abelian variety $A^{n+1}$.
Note that for $k\geq 2$ the zero locus of generic section of $L^{\otimes k}$ is smooth for any abelian variety by Bertini theorem. 
The cobordism class $[\Theta^n(k)]$ does not depend on the choice of such section.

\begin{Theorem} \cite{BV-2024} The cobordism class $[\Theta^n(k)]$ can be expressed as
\beq{cobt}
[\Theta^n(k)]=k^{n+1}[\Theta^{n}].
\eeq
\end{Theorem}

The proof follows from the identification of the normal bundle to $\Theta^n(k)$  with $\mathcal L^{\otimes k},$ where $\mathcal L=i^*(L)$.

As  a corollary we have the formulas for the generating functions
$$
C^\tau(\Theta^n(k),t)=k^{n+1} (n+1)! M_n(t_1,\dots,t_n), \, \, C^\nu(\Theta^n(k),t)=k^{n+1} (n+1)! t_n,
$$
so $\Theta^n(k)$ is also numerically extremely divisible with the Euler characteristic $\chi(\Theta^n(k))=k^{n+1}(n+1)!.$

Let us discuss in more detail what happens in two dimensions.
From Thom-Hirzebruch signature theorem
$$
c_1^2(X)=2\chi(X)+3\tau(X)
$$
it follows that in that case the Euler characteristic $\chi(X)$ must divide $3\tau(X)$, where $\tau(X)$ is the signature of $X.$

 In general, since in two dimensions the Chern numbers are some integer linear combinations of $c_1^2$ and $c_2$, we are looking for surfaces $X$ with $c_1^2(X)$ divisible by $c_2(X)=\chi(X)$.
 
For the smooth degree $d$ del Pezzo surfaces $S_d, d\leq 9$ we have the following
 
\begin{prop} 
All the Chern numbers of a del Pezzo surface $S_d$ are divisible by its Euler characteristic if and only if $d=6,8,9.$
\end{prop}

Indeed, $c_1^2(S_d)=d$ is divisible by $\chi(S_d)=12-d$ only in those cases. Note that $d=9$ corresponds to $\mathbb CP^2,$ while $d=8,6$ - to $\mathbb CP^2$ blown up at 1 and 3 points respectively.

For a smooth toric surface $X_N$ corresponding to $N$-gon $P_N$ we have
\begin{prop} 
All the Chern numbers of toric surface $X_N$ are divisible by its Euler characteristic if and only if $N=3,4,6,12.$
\end{prop}

Indeed, it is known that Euler characteristic of any smooth toric variety equals the number of vertices of the corresponding simple polytope, so in our case $c_2(X_N)=\chi(X_N)=N.$ 
The Todd genus in 2D case can be expressed as
$
Td(X)=\frac{1}{12}(c_1^2+c_2)
$
and for any toric variety is known to be 1.
Thus we have
$
c_1^2(X_N)=12-N,
$
which is divisible by $N$ if and only if $N$ divides 12.

Note that the cobordism class of the corresponding surface $X_N$ depends only on $N$, but not on the choice of the polygon $P_N.$
In particular, for $N=3$ and $N=4$ one can choose $\mathbb CP^2$ and $(\mathbb CP^1)^2$ respectively, for $N=6$ - permutohedral surface $X_{\Pi}$ (which is cobordant to $\Theta^2$, see \cite{BV-2022}).

Similarly, one can check that the degree $d$ hypersurfaces $V_d \subset \mathbb CP^3$ are numerically extremely divisible only for $d=1,2,4$ with $V_4$ being $K3$ surface.

Note that the $K3$ surfaces with $c_1^2=0, \, c_2=24$ and the Enriques surfaces with $c_1^2=0, \, c_2=12$, as well as Calabi-Yau 3-folds with $c_1^3=c_1c_2=0$, are numerically extremely divisible for obvious reasons.

For the compact complex surfaces of general type $X$  the {\it Bogomolov-Miyaoka-Yau inequality} 
$
c_1^2(X)\leq 3c_2(X)
$
implies that the {\it Chern slope} $$s=\frac{c_1^2(X)}{c_2(X)} \leq 3,$$ so in our case $s$ can be equal only to 1,2,3, with all these values realisable by some algebraic surface of general type (see Ch. VII in \cite{Barth}). Note that the theta divisor $\Theta^2$ with $c_1^2=c_2=6$ has the slope $s=1.$

 In particular, the surfaces $X$ with the extremal Chern invariants, satisfying the equality $c_1^2(X)=3c_2(X),$ are numerically extremely divisible in our sense. They include the fake projective planes \cite{PY}, which were actively studied after Mumford \cite{Mumford} (see recent results in \cite{BK}). Yau \cite{Yau} proved that  
$X$ must be isomorphic to a quotient of the unit ball in $\mathbb C^2$ by an infinite discrete group. 
More details about geography of the Chern numbers in dimensions 2 and 3 can be found in \cite{Barth, Hunt}.

In higher dimensional toric case we have the numerically extremely divisible examples of $\mathbb CP^n$ and $(\mathbb CP^1)^n$ corresponding to simplex and cube, but for the permutohedral variety $X_\Pi^3$ we have $c_{(3)}(X_\Pi^3)=20$, which is not divisible by $\chi (X_\Pi^3)=24.$ 
We have also the example of 3D flag variety $F_3=U(3)/T^3$ with the Chern numbers
$
c_3=\chi=6, \, c_1c_2=24, \, c_1^3=48
$
(see \cite{BT-2007}). 

How exceptional are such examples in higher dimensions is still to be understood.

It is also worthy to look at a more general class of irreducible algebraic varieties $M^n$ with $d(M^n)>1.$ The calculation and possible geometric meaning of $d(M^n)$ is another interesting question. 

Hirzebruch \cite{Hirzebruch_1958} studied a related question of describing pairs of the complex manifolds $M^n_1$ and $M^n_2$ having the proportional Chern numbers
$$
c_\lambda(M_1^n)=\mu \, c_\lambda(M_2^n), \,\, \mu \in \mathbb  Z
$$
for all partitions $\lambda$ of $n.$
The corresponding complex cobordism classes are proportional:
$
[M_1^n]=\mu \, [M_2^n].
$
Hirzebruch proved that this is the case for the compact quotients of a bounded symmetric domain and their compact duals  \cite{Hirzebruch_1958} (see also \cite{Borel_1963}).

We have already seen other examples of such pairs: $K3$ and Enriques surfaces with the relation $[K3]=2[E]$ and the divisors in abelian variety satisfying $[\Theta^n(k)]=k^{n+1}[\Theta^n]$, but in general such pairs need not to be numerically extremely divisible. 

Note that the divisibility of all Chern numbers of algebraic variety $X^n$ by $d>1$ does not imply the existence of algebraic variety $Y^n$ such that $[X^n]=d [Y^n].$
The simplest example is $\mathbb CP^n$ with all Chern numbers divisible by $n+1$, but $[\mathbb CP^n]\neq (n+1)[Y^n]$ for any complex manifold $Y^n$ since the Todd genus of $\mathbb CP^n$ equals 1.

\section{Concluding remark}
The example of Catalan formal group and the related Hirzebruch genus shows potential importance of the formal groups with the exponential being a polynomial, for example,
$$
f_k(z)=z-z^k.
$$
As we have seen, when $k=2$ the corresponding logarithm is a generating function of the Catalan numbers. For $k=3$ one can use the Lagrange inversion formula to show that we get a generating function of the Fuss-Catalan numbers $F_n$:
$$
g_3(u)=\sum_{n=0}^\infty F_n u^{2n+1}, \quad F_n=\frac{1}{2n+1}{3n \choose n}.
$$
We plan to discuss the topological aspects of such formal groups and other related questions in a separate publication.

\section{Acknowledgements.}
We are very grateful to Artie Prendergast for the helpful discussions and to Tom Copeland for his comments on the earlier version of this paper.

\end{document}